\documentclass[12pt,journal,final,onecolumn]{IEEEtran}

\makeatletter

\let\proof\@undefined
\let\endproof\@undefined
\makeatother

\usepackage[dvips]{graphicx}
\usepackage{epsfig}
\usepackage[cmex10]{amsmath}
\usepackage{amssymb}
\usepackage{amsfonts}
\usepackage{amsthm}
\usepackage{bm}
\usepackage{color}
\usepackage[mathscr]{eucal}
\usepackage{algorithm}
\usepackage[noend]{algorithmic}
\usepackage{times}
\usepackage{balance}

\newcommand{\mc}[1]{\mathcal{#1}}

\newcommand{\defeq}{\triangleq}
\newcommand{\Pp}{\mathbb{P}}

\newcommand{\E}{\mathbb{E}}

\newtheorem{lemma}{Lemma}

\newtheorem{theorem}[lemma]{Theorem}
\newtheorem{corollary}[lemma]{Corollary}

\theoremstyle{definition}
\newtheorem{exampledummy}{Example}
\newenvironment{example}{%
    \begin{exampledummy}%
        \upshape%
        
}
{\qed%
\end{exampledummy}}

\theoremstyle{remark}
\newtheorem*{remark}{Remark}

\title{Tracking Stopping Times \\Through Noisy Observations}

\author{Urs Niesen and Aslan Tchamkerten 
\thanks{The authors are with the Massachusetts Institute of Technology,
Cambridge, MA 02139.
Email: \{uniesen,tcham\}@mit.edu}
\thanks{This work was supported in part by NSF under Grant
No.~CCF-0515122, by DoD MURI Grant No. N00014-
07-1-0738, and by a University IR\&D Grant from Draper
Laboratory.}
}

\begin{document}

\maketitle

\begin{abstract} 
    A novel quickest detection setting is proposed which is a
    generalization of the well-known Bayesian change-point detection model. Suppose
    $\{(X_i,Y_i)\}_{i\geq 1}$ is a sequence of pairs of random
    variables, and that $S$ is a stopping time with respect to
    $\{X_i\}_{i\geq 1}$. The problem is to find a stopping time $T$ with
    respect to $\{Y_i\}_{i\geq 1}$ that optimally tracks $S$, in the
    sense that $T$ minimizes the expected {\emph{reaction delay}} $\E
    (T-S)^+$, while keeping the {\emph{false-alarm probability}}
    $\Pp(T<S)$ below a given threshold $\alpha \in [0,1]$. This problem
    formulation applies in several areas, such as in communication, detection,
    forecasting, and quality control.

    Our results relate to the situation where the $X_i$'s and $Y_i$'s  take values in finite
    alphabets and where $S$ is
    bounded by some positive integer $\kappa$. By using elementary
    methods based on the analysis of the tree structure of stopping
    times, we exhibit an algorithm that computes the optimal
    average reaction delays for all $\alpha \in [0,1]$, and constructs
    the associated optimal stopping times $T$. Under certain conditions on 
    $\{(X_i,Y_i)\}_{i\geq 1}$ and $S$, the algorithm running time is polynomial in $\kappa$.
\end{abstract}

{\keywords Algorithm, quickest detection problem, decision theory, synchronization, forecasting, monitoring, sequential analysis
}

\normalsize

\section{Problem Statement}
\label{sec:intro}

The tracking stopping time (TST) problem is defined as follows. Let
$\{(X_i,Y_i)\}_{i\geq 1}$ be a sequence of pairs of random variables.  Alice
observes $X_1,X_2,\ldots$ and  chooses a stopping time (s.t.) $S$ with respect
to that sequence.\footnote{ Recall that a stopping time with respect to a
sequence of random variables $\{X_i\}_{i\geq 1}$  is a random variable $S$
taking values in the positive integers such that the event $\{S=n\}$,
conditioned on $\{X_i\}_{i=1}^{n}$, is independent of
$\{X_{i}\}_{i=n+1}^{\infty}$ for all $n\geq 1$.  A  stopping time $S$ is
\emph{non-randomized} if $\Pp(S=n\vert X^n = x^n) \in \{0,1\}$ for all $x^n\in
\mc{Y}^n$ and $n\geq 1$.  A stopping time $S$ is \emph{randomized} if
$\Pp(S=n\vert X^n=x^n)\in [0,1]$ for all $x^n\in \mc{X}^n$ and $n\geq 1$. }
Knowing the distribution of $\{(X_i,Y_i)\}_{i\geq 1}$ and the stopping rule
$S$, but having access only to the $Y_i$'s, Bob wishes to find a s.t. that
gets as close as possible to Alice's. Specifically, Bob aims to find a s.t. $T$
minimizing the expected reaction delay $\E(T - S)^+ \defeq \E\max\{0,T-S\}$,
while keeping the false-alarm probability $\Pp(T<S)$ below a certain threshold
$\alpha\in[0,1]$.

\begin{example} {\bf Monitoring} 
    \label{eg:detection}

    Let $X_i$ be the distance of an object from a barrier at time $i$,
    and let $S$ be the first time the object hits the barrier, i.e.,
    $S\defeq\inf\{i\geq 1\,:\, X_i=0\}$. Assume we have access to $X_i$
    only through a noisy measurement $Y_i$,  and that we want to raise an
    alarm as soon as the  object hits the barrier.  This problem can be
    formulated as the one of finding a s.t. $T$ with respect to
    the $Y_i$'s that minimizes the expected reaction delay $\E(T-S)^+$,
    while keeping the false-alarm probability $\Pp(T<S)$ small enough. 
\end{example}

Another situation where the TST problem applies is in the context
of communication over channels with feedback.  Most of the studies
related to feedback communication assume perfect feedback, i.e., the
transmitter is fully aware of the output of the channel as observed by
the receiver. Without this assumption --- i.e., if the feedback link is
noisy --- a synchronization problem may arise between the transmitter and
the receiver which can be formulated as a TST problem, as shown in the following
example. 
\begin{example} {\bf Communication}
    \label{eg:communication}

    It is well known that the presence of a noiseless feedback link
    allows to dramatically increase the reliability for a given
    communication delay (see, e.g.,~\cite{H}). However, to take
    advantage of feedback, variable length codes are often
    necessary.\footnote{The reliability function associated with block
    coding schemes is lower than the one associated with variable length
    coding. For symmetric channels, for instance, the reliability
    function associated with block coding schemes is limited by the
    sphere packing bound, which is lower than the best optimal error
    exponent attainable with variable length coding (\cite{B,D}). } This
    can be observed by looking at a non-perfect binary erasure channel.
    In this case, any block coding strategy yields a strictly positive
    error probability. In contrast, consider the variable length
    strategy where the encoder keeps sending the bit it wishes to convey
    until it is successfully received. This simple strategy achieves
    error-free communication at a rate equal to the capacity of the
    channel in question.  Can we still use this coding strategy if the
    feedback channel is (somewhat) noisy?  Because of the noisy feedback
    link, a synchronization problem between the decoder and the encoder
    arises: if the first non-erased output symbol occurs at time $S$,
    what should be sent at time $S+1$? This agreement problem occurs
    because the encoder observes now only a noisy version of the symbols
    received by the decoder (see Fig.~\ref{fig:becbec}). In particular, the
    first non-erased output symbol may not be recognized as such by the
    encoder.\footnote{For fixed length coding strategies over channels
    with noisy feedback we refer the reader to \cite{KLW,BY}.} 

    \begin{figure}
        \begin{center}
            \input{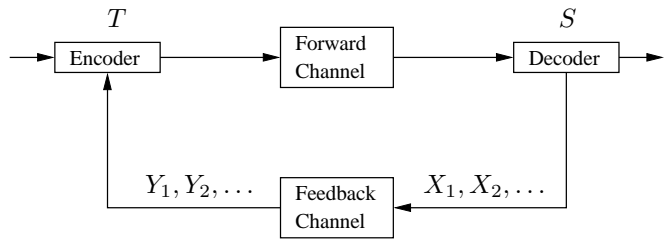}
        \end{center}
        \caption{
        \label{fig:becbec} 
        The decoding time $S$ depends on the output of the forward channel. The
        encoder decides to stop transmission at time $T$ based on the output of
        the feedback channel. If the feedback channel is noisy, $S$ and $T$ need
        not coincide.
        }

    \end{figure} 

    Instead of treating the synchronization issue that results from the
    use of variable length codes over channels with noisy feedback, let
    us consider the simpler problem of finding the minimum delay needed
    by the encoder to realize that the decoder has made a decision. In
    terms of the TST problem, Alice and Bob represent the decoder and
    the encoder, the $X_i$'s and $Y_i$'s correspond to the input and
    output symbols of the feedback channel, whereas $S$ and $T$
    represent the decoding time and the time the encoder stops
    transmission, respectively. Here $\E(T-S)^+$ represents the delay
    needed by the encoder to realize that the decoder has made a
    decision, and we aim to minimize it given that the probability of
    stopping transmission too early, $\Pp(T<S)$, is kept below a certain
    threshold $\alpha$. 

    Note that, in the context of feedback communication, it would be
    reasonable to define the communication rate with respect to the
    overall delay $S+(T-S)^+=\max\{S,T\}$. This definition, in contrast with the one
    that takes into account only the decoding time (such as for rateless codes),
    puts the delay constraint on both the transmitter and the receiver. In the
    Example \ref{eg:continu}, we investigate the highest achievable rate with
    respect to the overall communication delay if the ``send until
    a non-erasure occurs'' strategy is used and both the forward and the
    feedback channels are binary erasure.
\end{example}

\begin{example} {\bf Forecasting}
    \label{eg:forecasting}

    A large manufacturing machine breaks down as soon as its cumulative
    fatigue hits a certain threshold. Knowing that a machine replacement
    takes, say, ten days, the objective is to order a new machine so that it
    is operational at the time the old machine breaks down. This prevents
    losses due to an interrupted manufacturing process as well as storage
    costs caused by an unused backup machine.

    The problem of determining the operating start date of the new machine
    can be formulated as follows. Let $X_n$ be the cumulative fatigue up to
    day $n$ of the current machine, and let $S$ denote the first day $n$
    that $X_n$ crosses the critical fatigue threshold. Since the replacement
    period is ten days, the first day $T$ a new machine is operational
    can be scheduled only on the basis of a (possibly randomized) function
    of $\{X_i\}_{i=1}^{T-10}$. By defining $Y_i$ to be equal to $X_{i-10}$
    if $i>10$ and else equal to zero, the day $T$ is now a s.t. with respect
    to $\{Y_i\}_{i\geq 1}$, and  we can formulate the requirement on $T$ as
    aiming to minimize $\E(T-S)^+$ while keeping $\Pp(T<S)$ below a certain
    threshold.  
\end{example}

Note that, in the forecasting example, in contrast with the monitoring
and communication examples, Alice has access to more information than
Bob. From the process she observes, she can deduce Bob's observations
--- here simply by delaying hers. This feature may be interesting in other
applications. The general formulation where Alice has access to more
information than Bob is obtained by letting the observation available to
Alice at time $i$ be $X_i=(\tilde{X}_i,\tilde{Y}_i)$, and the
observation available to Bob be $Y_i = \tilde{Y}_i$.
\begin{example} {\bf Bayesian Change-Point Detection} 
    \label{eg:cp}

    In this Example we will see how the TST setting generalizes the Bayesian
    version of the change-point detection problem, a long studied problem with
    applications to industrial quality control and that dates back to the
    $1940$'s \cite{AGP}. The Bayesian change-point problem is formulated as
    follows.  Let $\theta$ be a random variable taking values in the positive
    integers.  Let $\{Y_i\}_{i\geq 1}$ be a sequence of random variables such
    that, given the value of $\theta$, the conditional probability of $Y_n$
    given $Y^{n-1}\defeq \{Y_i\}_{i=1}^{n-1}$ is $P_0(\cdot|Y^{n-1})$ for $n <
    \theta$ and is $P_1(\cdot|Y^{n-1})$ for $n\geq \theta$.  We are interested
    in a s.t.  $T$ with respect to the $Y_i$'s minimizing the change-point
    reaction delay $\E(T-\theta)^+$, while keeping the false-alarm probability
    $\Pp(T<\theta)$ below a certain threshold $\alpha\in[0,1]$. 

    Shiryaev~(see, e.g.,\cite{S2},\cite[Chapter 4.3]{S}) considered the Lagrangian
    formulation of the above problem: Given a constant $\lambda\geq
    0$, minimize  $$\E(T-\theta)^++\lambda \Pp(T<\theta)$$ over all s.t.'s $T$.
    Assuming a geometric prior on the change-point $\theta$ and that
    before and after $\theta$ the observations are independent with
    common density function $f_0$ for $t<\theta$ and $f_1$ for $t\geq
    \theta$, Shiryaev showed that the optimal $T$ stops as soon as the
    posterior probability that a change occurred exceeds a certain fixed
    threshold.  Later Yakir~\cite{Y} generalized Shiryaev's result by
    considering finite-state Markov chains. For more general prior
    distributions on $\theta$, the problem is known to become difficult
    to handle. However, in the limit $\alpha \rightarrow 0$,
    Lai~\cite{L} and, later, Tartakovsky and Veeravalli~\cite{TV},
    derived asymptotically optimal detection policies for the Bayesian
    change-point problem  under general assumptions on the distributions
    of the change-point and observed process.\footnote{For the
    non-Bayesian version of the change-point problem we refer the reader
    to \cite{Lo,Po,Mo}.}

    To see that the Bayesian change-point problem can be formulated as a TST
    problem, it suffices to define the sequence of binary random variables
    $\{X_i\}_{i\geq 1}$ such that $X_i=0$ if $i<\theta$ and $X_i=1$ if
    $i\geq\theta$, and to let $S\defeq\inf\{i:X_i=1\}$ (i.e., $S=\theta)$. The
    change-point problem defined by $\theta$ and $\{Y_i\}_{i\geq 1}$ becomes the TST
    problem defined by $S$ and $\{(X_i,Y_i)\}_{i\geq 1}$.
    However, the TST problem cannot, in general,  be formulated as a Bayesian
    change-point problem. Indeed, the Bayesian change-point problem yields
    for any $k>n$ 
    \begin{align}
        \label{eq:cpdiff}
        \Pp(\theta=k\vert & Y^n=y^n, \theta>n) \nonumber \\
        & = \frac{\Pp(Y^n=y^n, \theta>n\vert \theta=k)\Pp(\theta=k)}
        {\Pp(Y^n=y^n\vert \theta>n)\Pp(\theta>n)} \nonumber\\
        & = \frac{\Pp(Y^n=y^n\vert \theta=k)\Pp(\theta=k)}
        {\Pp(Y^n=y^n\vert \theta>n)\Pp(\theta>n)} \nonumber\\
        & = \Pp(\theta=k\vert \theta>n)
    \end{align}
    since $\Pp(Y^n=y^n\vert \theta=k)=\Pp(Y^n=y^n\vert \theta>n)$.  Therefore,
    conditioned on the event $\{\theta>n\}$, the first $n$ observations
    $Y^n$ are independent of $\theta$. In other words, given that no change
    occurred up to time $n$, the observations $y^n$ are useless in
    predicting the value of the change-point $\theta$. In contrast, for the
    TST problem, in general we have 
    \begin{align}
        \Pp(S=k\vert Y^n=y^n, S>n)
        & \ne \Pp(S=k\vert S>n)
    \end{align}
    because $\Pp(Y^n=y^n\vert S=k)\ne \Pp(Y^n=y^n\vert S>n)$.
\end{example}

As is argued in the last example, the TST problem is a generalization of
the Bayesian change-point problem, which itself is analytically
tractable only in special cases. This makes an analytical treatment of
the general TST problem difficult. Instead, we present an algorithmic
solution to this problem for an arbitrary process $\{(X_i,Y_i)\}_{i\geq
1}$ and an arbitrary stopping time $S$ bounded by some constant
$\kappa\geq 1$.  The proof of correctness of this algorithm provides
insights into the structure of the optimal stopping time $T$ tracking
$S$, and into the tradeoff between expected delay $\E(T-S)^+$ and
probability of false-alarm $\Pp(T<S)$. Under some conditions on
$\{(X_i,Y_i)\}_{i\geq 1}$ and $S$, the computational complexity of this
algorithm is polynomial in $\kappa$.

The rest of the paper is organized as follows.  In Section
\ref{sec:det_stop_tim}, we provide some basic properties of the TST
problem defined over a finite alphabet process $\{(X_i,Y_i)\}_{i\geq
1}$, and in Section \ref{sec:mainresult} we provide an algorithmic
solution to it.  In Section~\ref{sec:perm}, we derive conditions under
which the algorithm has low complexity and illustrate this in
Section~\ref{sec:lookahead} with examples.

\section{The Optimization Problem}
\label{sec:det_stop_tim}

Let $\{(X_i,Y_i)\}_{i\geq 1}$ be a discrete-time process where the
$X_i$'s and $Y_i$'s take value in some finite alphabets $\cal{X}$ and
$\cal{Y}$, respectively. Let $S$ be a s.t.  with respect to
$\{X_i\}_{i\geq 1}$ such that $S\leq \kappa$ almost surely for some constant
$\kappa \geq 1$. We aim to find for any $\alpha \in [0,1]$ 
\begin{equation} 
    \label{eq:prob}
    d(\alpha) \defeq \min_{\substack{T:\Pp(T<S)\leq\alpha\\ T\leq \kappa}}\E(
    T-S)^{+} 
\end{equation} 
where the s.t.'s $T$ are possibly randomized. Note that the restriction $T\leq
\kappa$ induces no loss of optimality.

\begin{figure}
    \begin{center}
        \input{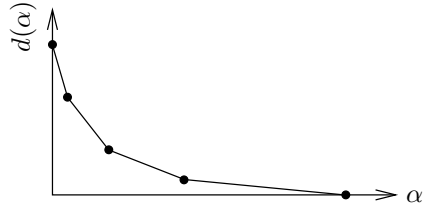}
    \end{center}
    \caption{Typical shape of the expected delay $d(\alpha)$ as a function
    of false-alarm probability $\alpha$. The break-points are achieved by
    non-randomized stopping times.}
    \label{fig:betalpha}
\end{figure}

Now, the set of all s.t.'s over $\{Y_i\}_{i\geq 1}$ is convex, and its
extreme points are non-randomized s.t.'s (\cite{BC},~\cite{EMS}). This
implies that any randomized s.t.  $T\leq\kappa$ can be written as a
convex combination of non-randomized s.t.'s bounded by $\kappa$, i.e.
\begin{equation*}
    \Pp(T=k)=\sum_j w_j\Pp(T_j=k)
\end{equation*}
for any integer $k$, where $\{T_j\}$ denotes the finite set of all
non-randomized s.t.'s bounded by $\kappa$, and where the $w_j$'s are nonnegative
and sum to one. Hence, because false-alarm and expected reaction delay can be written as
\begin{align*}
    \Pp(T<S) & = \sum_j w_j \Pp(T_j<S) \\
    \E(T-S)^+ & = \sum_j w_j \E(T_j-S)^+\;,
\end{align*}
 the function $d(\alpha)$ is
convex and piecewise linear, with break-points achieved by
non-randomized s.t.'s. Its typical shape is depicted in
Figure~\ref{fig:betalpha}.

For $\lambda\geq 0$, define the Lagrangian
\begin{equation}
    \label{eq:lagrange}
    J_\lambda(T) \defeq \E(T-S)^++\lambda \Pp(T<S).
\end{equation}
\begin{lemma}
    We have
    \label{thm:dual}
    \begin{equation*}
        d(\alpha) = \sup_{\lambda\geq 0} \min_{T\leq \kappa} \left(
        J_{\lambda}(T)-\lambda
        \alpha \right),
    \end{equation*}
    where the minimization is over all non-randomized s.t.'s bounded by
    $\kappa$.
\end{lemma}
\begin{proof}
    The convex minimization problem in~\eqref{eq:prob} admits at least one feasible
    point, namely $T=\kappa$. Therefore strong Lagrange duality holds (see,
    e.g.,~\cite[Chapter 5]{BV}), and we obtain
    \begin{equation}
        \label{eq:dual}
        d(\alpha) = \sup_{\lambda\geq 0} \min_{T\leq \kappa} \left( J_{\lambda}(T)-\lambda
        \alpha \right).
    \end{equation} 
    Because $d(\alpha)$ is convex with extreme points achieved
    by non-randomized s.t.'s, we may restrict the minimization
    in~\eqref{eq:dual} to be over the set of non-randomized s.t.'s
    bounded by~$\kappa$. 
\end{proof}

\section{An Algorithm for Computing $d(\alpha)$}
\label{sec:mainresult}

We first establish a few preliminary results later used to evaluate
$\min_{T}J_\lambda(T)$. Emphasis is put on the finite tree representation of
bounded s.t.'s with respect to finite alphabet processes. We then
provide an algorithm that computes the entire curve $d(\alpha)$.

We introduce a few notational conventions. The set $\mc{Y}^*$
represents all finite sequences over $\cal{Y}$. An element in
$\mc{Y}^*$ is denoted either by $y^n$ or by $\bm{y}$, depending on
whether or not we want to emphasize its length. To any non-randomized
s.t. $T$, we associate a unique  $\vert\mc{Y}\vert$-ary tree $\mc{T}$
(i.e., all the nodes of $\mc{T}$ have either zero or exactly
$\vert\mc{Y}\vert$ children) having each node specified by some
$\bm{y}\in \mc{Y}^*$, where $\rho \bm{y}$ represents the vertex path
from the root $\rho$ to the node $\bm{y}$. The depth of a node
$y^n\in\mc{T}$ is denoted by $l(y^n) \defeq n$. The tree consisting only
of the root is the trivial tree.  A node $y^n\in \mc{T}$ is a leaf if
$\Pp(T=n\vert Y^n = y^n) = 1$. We denote by $\mc{L}(\mc{T})$ the leaves
of $\mc{T}$ and by $\mc{I}(\mc{T})$ the intermediate (or non-terminal)
nodes of $\mc{T}$.  The notation $T(\mc{T})$ is used to denote the
(non-randomized) s.t. $T$ induced by the tree $\mc{T}$.  Given a node
$\bm{y}$ in $\mc{T}$, let $\mc{T}_{\bm{y}}$ be the subtree of $\mc{T}$
rooted in $\bm{y}$. Finally, let $\mc{D}(\mc{T}_{\bm{y}})$ denote the
descendants of $\bm{y}$ in $\mc{T}$. The next example illustrates these
notations.

\begin{example}
    \label{eg:tree}
    Let $\mc{Y} = \{0,1\}$ and $\kappa = 2$. The tree $\mc{T}$ depicted in
    Figure~\ref{fig:tree} corresponds to the non-randomized s.t.
    $T$ taking value one if $Y_1 = 1$ and value $2$ if $Y_1 = 0$. The
    sets $\mc{L}(\mc{T})$ and $\mc{I}(\mc{T})$ are given by
    $\{00,01,1\}$ and $\{\rho, 0\}$, respectively. The subtree
    $\mc{T}_0$ of $\mc{T}$ consists of the nodes $\{0,00,01\}$, and its
    descendants $\mc{D}(\mc{T}_{0})$ are $\{00,01\}$. The
    subtree $\mc{T}_{\rho}$ is the same as $\mc{T}$, and its descendants
    $\mc{D}(\mc{T}_{\rho})$ are $\{0,1,00,01\}$.
    \begin{figure}[h!]
        \begin{center}
            \input{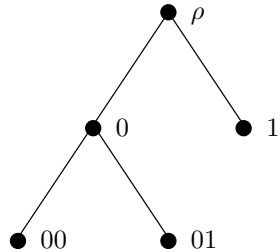}
        \end{center}
        \caption{Tree corresponding to the s.t. $T$ defined by $T=1$ if
        $Y_1=1$, and $T=2$ else.}
        \label{fig:tree}
    \end{figure}
\end{example}

Below, we describe an algorithm that, for a given s.t. $S$, constructs a
sequence of s.t.'s $\{T(\mc{T}^m)\}_{m=0}^{M}$ and Lagrange multipliers
$\{\lambda_m\}_{m=0}^{M}$ with the following two properties.  First, the
$\mc{T}^m$'s and $\lambda_m$'s are ordered in the sense that
$\mc{T}^{M}\subset\mc{T}^{M-1}\subset\ldots\subset\mc{T}^{0}$ and
$0=\lambda_{M}\leq\lambda_{M-1}\leq\ldots\leq\lambda_1\leq\lambda_0=\infty$.
(Here the symbol $\subset$ denotes inclusion, not necessarily strict.)
Second, for any $m\in\{0,\ldots,M\}$ and $\lambda\in
(\lambda_{m},\lambda_{m-1}]$ the tree $\mc{T}^{m-1}$ minimizes
$J_{\lambda}(\mc{T})\defeq J_{\lambda}(T(\mc{T}))$ among all
non-randomized s.t.'s. The algorithm builds upon ideas from the CART
algorithm for the construction of classification and regression
trees~\cite{BFO}. 

Before we state the algorithm, we need to introduce a few quantities.
Given a non-randomized s.t. $T$ represented by its
$\vert\mc{Y}\vert$-ary tree $\mc{T}$, we write the Lagrangian
$J_{\lambda}(\mc{T})$ as
\begin{align*}
    J_{\lambda}(\mc{T}) 
    & = \E(T-S)^++\lambda\Pp(T<S) \\
    & = \sum_{\bm{y}\in\mc{L}(\mc{T})}\Pp(\bm{Y}=\bm{y})
    \Big(\E\big( (l(\bm{y})-S)^+ \vert\bm{Y}=\bm{y}\big) \\
    & \quad +\lambda\Pp\big(S>l(\bm{y})\vert\bm{Y}=\bm{y}\big)\Big) \\
    & = \sum_{\bm{y}\in\mc{L}(\mc{T})} b(\bm{y})+\lambda a(\bm{y}) \\
    & = \sum_{\bm{y}\in\mc{L}(\mc{T})} J_{\lambda}(\bm{y}),
\end{align*}
where 
\begin{align*}
    a(\bm{y})&\defeq
    \Pp(\bm{Y}=\bm{y}) \Pp(S>l(\bm{y})\vert\bm{Y}=\bm{y}),  \\
    b(\bm{y})&\defeq \Pp(\bm{Y}=\bm{y}) \E\big( (
    l(\bm{y})-S)^+ \vert\bm{Y}=\bm{y}\big), \\
    J_{\lambda}(\bm{y})&\defeq b(\bm{y})+\lambda a(\bm{y})\;.
\end{align*}
We extend the definition of $J_{\lambda}(\cdot)$ to subtrees of $\mc{T}$
by setting 
\begin{equation*}
    J_{\lambda}(\mc{T}_{\bm{y}}) \defeq
    \sum_{\bm{\gamma}\in\mc{L}(\mc{T}_{\bm{y}})}J_{\lambda}(\bm{\gamma}).
\end{equation*}
With this definition\footnote{We used $T$, $\mc{T}$, $\mc{T}_{\bm{y}}$,
and $\bm{y}$, as possible arguments of $J_{\lambda}(\cdot)$. No
confusion should arise from this slight abuse of notation, since for
non-randomized s.t.'s all of these arguments can be interpreted as trees.}
\begin{equation*}
    J_{\lambda}(\mc{T}_{\bm{y}}) = 
    \left\{ \begin{array}{ll}
        J_{\lambda}(\bm{y}) & \mbox{if $\bm{y}\in\mc{L}(\mc{T})$},\\
        \sum_{\gamma\in\mc{Y}}J_{\lambda}(\mc{T}_{\bm{y}\gamma}) 
        & \mbox{if $\bm{y}\in\mc{I}(\mc{T})$}.
    \end{array} \right. 
\end{equation*}
Similarly, we define
\begin{align*}
    a(\mc{T}_{\bm{y}}) & \defeq \sum_{\bm{\gamma}\in \mc{L}(\mc{T}_{\bm{y}})}a(\bm{\gamma}), \\
    b(\mc{T}_{\bm{y}}) & \defeq \sum_{\bm{\gamma}\in \mc{L}(\mc{T}_{\bm{y}})}b(\bm{\gamma}).
\end{align*}

For a given $\lambda\geq 0$ and $\mc{T}$, define
$\mc{T}(\lambda)\subset\mc{T}$ to be the subtree of $\mc{T}$ having the
same root, and such that $J_{\lambda}(\mc{T}(\lambda))\leq
J_{\lambda}(\mc{T}')$ for all subtrees (with same root) $\mc{T}'\subset
\mc{T}$, and $\mc{T}(\lambda) \subset \mc{T}'$ for all subtrees (with
same root) $\mc{T}'\subset \mc{T}$ satisfying
$J_{\lambda}(\mc{T}(\lambda))= J_{\lambda}(\mc{T}')$.  In words, among
all subtrees of $\mc{T}$ yielding a minimal cost for a given $\lambda$,
the tree $\mc{T}(\lambda)$ is the smallest. As we shall see in
Lemma~\ref{thm:optimal}, such a smallest subtree always exists, and
hence $\mc{T}(\lambda)$ is well defined.

\begin{remark}
    Note that $\mc{T}_{\bm{y}}(\lambda)$ is different from
    $(\mc{T}(\lambda))_{\bm{y}}$. Indeed, $\mc{T}_{\bm{y}}(\lambda)$ refers
    to the optimal subtree of $\mc{T}_{\bm{y}}$ with respect to $\lambda$,
    whereas $(\mc{T}(\lambda))_{\bm{y}}$ refers to subtree rooted in
    $\bm{y}$ of the optimal tree $\mc{T}(\lambda)$.
\end{remark}

\begin{example}
    Consider again the tree $\mc{T}$ in Figure~\ref{fig:tree}. Assume
    $J_\lambda(\rho)=4, J_\lambda(0)=2, J_\lambda(1)= J_\lambda(00)= J_\lambda(01)=1$.
    Then 
    \begin{align*}
        J_\lambda(\mc{T}) & = J_\lambda(1)+J_\lambda(00)+ J_\lambda(01)=3,\\
        J_\lambda(\mc{T}_0) & = J_\lambda(00)+ J_\lambda(01)=2.
    \end{align*}
    The smallest optimal subtree of $\mc{T}$ having the same 
    root is $\mc{T}(\lambda) = \{\rho,0,1\}$ and
    \begin{equation*}
        J_\lambda(\mc{T}(\lambda)) = J_\lambda(0)+J_\lambda(1) = 3.
    \end{equation*}
    The smallest optimal subtree of $\mc{T}_0$ having the same
    root is $\mc{T}_0(\lambda) = \{0\}$ and
    \begin{equation*}
        J_\lambda(\mc{T}_0(\lambda)) = J_\lambda(0) = 2.
    \end{equation*}
\end{example}

Given a $\vert\mc{Y}\vert$-ary tree $\mc{T}$, and $\lambda\geq 0$, the
following lemma shows that $\mc{T}(\lambda)$ always exists and
characterizes $\mc{T}(\lambda)$ and $J_\lambda(\mc{T}(\lambda))$. The
reader may recognize the finite-horizon backward induction algorithm whose
detailed proof can be found in standard textbooks (e.g.,~\cite[Chapter 3
and 4]{CRS}). 

\begin{lemma}
    \label{thm:optimal}
    Given a $\vert\mc{Y}\vert$-ary tree $\mc{T}$ and $\lambda\geq 0$.
    For every $\bm{y}\in\mc{I}(\mc{T})$,
    \begin{equation*}
        J_{\lambda}(\mc{T}_{\bm{y}}(\lambda)) = \min\{J_{\lambda}(\bm{y}),
        \sum_{\gamma\in\mc{Y}}J_{\lambda}(\mc{T}_{\bm{y}\gamma}(\lambda))\},
    \end{equation*}
    and 
    \begin{equation*}
        \mc{T}_{\bm{y}}(\lambda) = 
        \begin{cases}
            \{\bm{y}\} 
            & \text{\!\!if $J_{\lambda}(\bm{y}) \leq 
            \sum_{\gamma\in\mc{Y}}J_{\lambda}(\mc{T}_{\bm{y}\gamma}(\lambda))$} \\
            \{\bm{y}\}\cup_{\gamma\in\mc{Y}}\mc{T}_{\bm{y}\gamma}(\lambda)
            & \text{\!\!else}.
        \end{cases}
    \end{equation*} 
    The optimal tree $\mc{T}(\lambda)$ and the corresponding cost
    $J_\lambda(\mc{T}(\lambda))$ are given by
    $J_{\lambda}(\mc{T}_{\bm{y}}(\lambda))$ and $\mc{T}_{\bm{y}}(\lambda)$ evaluated at
    $\bm{y}=\rho$.
\end{lemma}
\begin{proof}
    By induction on the depth of the tree starting from the root. 
\end{proof}
From the structure of the cost function $J_\lambda(\cdot)$, the larger
the value of $\lambda$, the higher the penalty on the error probability.
Therefore one expects that the larger the $\lambda$ the ``later'' the
optimal tree $\mc{T}(\lambda)$ will stop.  Indeed, Lemma~\ref{thm:subtree}
states that the tree corresponding to the optimal s.t. of a
smaller $\lambda$ is a subtree of the tree corresponding to the optimal
s.t. of a larger $\lambda$.  In other words, if $\lambda\leq
\tilde{\lambda}$, in order to find $\mc{T}({\lambda})$, we can restrict
our search to subtrees of $\mc{T}(\tilde{\lambda})$.
\begin{lemma}
    \label{thm:subtree}
    Given a tree $\mc{T}$, if $\lambda\leq\tilde{\lambda}$ then
    $\mc{T}(\lambda)\subset\mc{T}(\tilde{\lambda})$.
\end{lemma}
\begin{proof}
    We have 
    \begin{align}
        \label{eq:alpharec}
        a(\mc{T}_{\bm{y}})
        & = \sum_{\bm{y}\bm{\gamma}\in\mc{L}(\mc{T}_{\bm{y}})} 
        \!\!\!\!\Pp(S>l(\bm{y}\bm{\gamma}) \vert Y^{l(\bm{y}\bm{\gamma})}=\bm{y}\bm{\gamma}) 
        \Pp(Y^{l(\bm{y}\bm{\gamma})}=\bm{y}\bm{\gamma}) \nonumber\\
        & \leq \sum_{\bm{y}\bm{\gamma}\in\mc{L}(\mc{T}_{\bm{y}})} 
        \!\!\!\!\Pp(S>l(\bm{y}) \vert Y^{l(\bm{y}\bm{\gamma})}=\bm{y}\bm{\gamma}) 
        \Pp(Y^{l(\bm{y}\bm{\gamma})}=\bm{y}\bm{\gamma}) \nonumber\\
        & = a(\bm{y}).
    \end{align}
    Similarly one shows that $b(\mc{T}_{\bm{y}}) \geq b(\bm{y})$. 

    By contradiction, assume $\lambda\leq\tilde{\lambda}$, but
    $\mc{T}(\lambda)$ is not a subset of $\mc{T}(\tilde{\lambda})$. Then
    there exists $\bm{y}\in\mc{L}(\mc{T}(\tilde{\lambda}))$ such that
    $\bm{y}\in\mc{I}(\mc{T}(\lambda))$. By definition of
    $\mc{T}(\tilde{\lambda})$ and Lemma \ref{thm:optimal}
    \begin{equation*}
        J_{\tilde{\lambda}}(\bm{y}) 
        \leq J_{\tilde{\lambda}}(\mc{T}_{\bm{y}}(\lambda)),
    \end{equation*}
    and thus
    \begin{equation}
        \label{eq:subtree1}
        b(\bm{y})+\tilde{\lambda}a(\bm{y})
        \leq b(\mc{T}_{\bm{y}}(\lambda))
        +\tilde{\lambda}a(\mc{T}_{\bm{y}}(\lambda)).
    \end{equation}
    Now, since $a(\mc{T}_{\bm{y}}(\lambda)) \leq a(\bm{y})$, 
    and $\lambda\leq\tilde{\lambda}$, 
    \begin{equation}
        \label{eq:subtree2}
        (\lambda-\tilde{\lambda})a(\bm{y})
        \leq (\lambda-\tilde{\lambda})a(\mc{T}_{\bm{y}}(\lambda)).
    \end{equation}
    Combining~\eqref{eq:subtree1} and~\eqref{eq:subtree2} yields
    \begin{equation*}
        b(\bm{y})+\lambda a(\bm{y})
        \leq b(\mc{T}_{\bm{y}}(\lambda))
        +\lambda a(\mc{T}_{\bm{y}}(\lambda)),
    \end{equation*}
    and therefore
    \begin{equation*}
        J_\lambda(\bm{y}) \leq J_\lambda(\mc{T}_{\bm{y}}(\lambda)).
    \end{equation*}
    Since $\bm{y}\in\mc{I}(\mc{T}(\lambda))$, this contradicts the
    definition of $\mc{T}(\lambda)$ by Lemma~\ref{thm:optimal}.
\end{proof}

The next theorem represents a key result.  Given a tree $\mc{T}$, it
characterizes the smallest value $\lambda$ can take for which
$\mc{T}(\lambda)=\mc{T}$. For a non-trivial tree $\mc{T}$, 
define for any $\bm{y}\in \mc{I}(\mc{T})$ 
\begin{equation*}
    g(\bm{y},\mc{T}) \defeq 
    \frac{b(\mc{T}_{\bm{y}})-b(\bm{y})}{a(\bm{y})-a(\mc{T}_{\bm{y}})},
\end{equation*}
where we set $0/0\defeq 0$.
The quantity $g(\bm{y},\mc{T})$ captures the tradeoff between the
reduction in delay $b(\mc{T}_{\bm{y}})-b(\bm{y})$ and the increase in
probability of false-alarm $a(\bm{y})-a(\mc{T}_{\bm{y}})$ if we stop at
some intermediate node $\bm{y}$ instead of stopping at the leaves
$\mc{L}(\mc{T}_{\bm{y}})$ of $\mc{T}$.

\begin{theorem}
    \label{thm:lambda}
    For any non-trivial tree $\mc{T}$ 
    \begin{align*}
        \inf\big\{\lambda\geq 0 : \mc{T}(\lambda) = \mc{T}\:\big\}
        =\max_{\bm{y}\in\mc{I}(\mc{T})} g\big(\bm{y},\mc{T}\big)\;.
    \end{align*}
\end{theorem}

\begin{proof}
    Let $\mc{T}$ be a non-trivial tree and $\bm{y}\in\mc{I}(\mc{T})$. We
    have
    \begin{align*}
        g(\bm{y},\mc{T}) 
        & = \frac{J_{\lambda}(\mc{T}_{\bm{y}})-\lambda a(\mc{T}_{\bm{y}})
        -J_{\lambda}(\bm{y})+\lambda a(\bm{y})}
        {a(\bm{y})-a(\mc{T}_{\bm{y}})} \\
        & = \frac{J_{\lambda}(\mc{T}_{\bm{y}})-J_{\lambda}(\bm{y})}
        {a(\bm{y})-a(\mc{T}_{\bm{y}})}+\lambda.
    \end{align*}
    By~\eqref{eq:alpharec}, $a(\mc{T}_{\bm{y}})\leq a(\bm{y})$, and hence
    the following implications hold:
    \begin{equation}
        \label{eq:gequiv}
        \begin{aligned}
            g(\bm{y},\mc{T}) \leq \lambda & \Longleftrightarrow 
            J_{\lambda}(\bm{y})\geq J_{\lambda}(\mc{T}_{\bm{y}}), \\
            g(\bm{y},\mc{T}) < \lambda & \Longleftrightarrow 
            J_{\lambda}(\bm{y}) > J_{\lambda}(\mc{T}_{\bm{y}}).
        \end{aligned}
    \end{equation}
    Therefore, if $\max_{\bm{y}\in\mc{I}(\mc{T})}g(\bm{y},\mc{T}) <
    \lambda  $  then
    \begin{equation}
        \label{eq:ineqj}
        J_{\lambda}(\bm{y})>J_{\lambda}(\mc{T}_{\bm{y}})
    \end{equation}
    for all $\bm{y}\in\mc{I}(\mc{T})$.

    We first show by induction that if
    \begin{equation*}
        \max_{\bm{y}\in\mc{I}(\mc{T})}g(\bm{y},\mc{T})<\lambda
    \end{equation*}
    then $\mc{T}(\lambda)=\mc{T}$. Consider a subtree of $\mc{T}$ having
    depth one and rooted in $\bm{y}$, say. Since by~\eqref{eq:ineqj} 
    $J_{\lambda}(\bm{y}) > J_{\lambda}(\mc{T}_{\bm{y}})$, 
    we have $\mc{T}_{\bm{y}}(\lambda) = \mc{T}_{\bm{y}}$ 
    by Lemma~\ref{thm:optimal}.  Now consider a subtree of
    $\mc{T}$ with depth $k$, rooted in a different $\bm{y}$, and assume the
    assertion to be true for all subtrees of $\mc{T}$ with depth up to
    $k-1$. In order to find $\mc{T}_{\bm{y}}(\lambda)$, we use
    Lemma~\ref{thm:optimal} and compare $J_{\lambda}(\bm{y})$ with
    $\sum_{\gamma\in\mc{Y}} J_{\lambda}(\mc{T}_{\bm{y}\gamma}(\lambda))$.
    Since $\mc{T}_{\bm{y}\gamma}$ is a subtree of $\mc{T}$ with depth less
    than $k$, we have $\mc{T}_{\bm{y}\gamma}(\lambda) =
    \mc{T}_{\bm{y}\gamma}$ by the induction hypothesis. Therefore
    \begin{equation*}
        \sum_{\gamma\in\mc{Y}} J_{\lambda}(\mc{T}_{\bm{y}\gamma}(\lambda))
        = \sum_{\gamma\in\mc{Y}} J_{\lambda}(\mc{T}_{\bm{y}\gamma})
        = J_{\lambda}(\mc{T}_{\bm{y}}),
    \end{equation*}
    and since $  J_{\lambda}(\mc{T}_{\bm{y}})<J_{\lambda}(\bm{y})$ by
    \eqref{eq:ineqj},  we have $\mc{T}_{\bm{y}}(\lambda)=\mc{T}_{\bm{y}}$ by
    Lemma~\ref{thm:optimal}, which concludes the induction step. Hence we
    proved that if $\max_{\bm{y}\in\mc{I}(\mc{T})}g(\bm{y},\mc{T}) <\lambda$, 
    then $\mc{T}(\lambda)=\mc{T}$.

    Second, suppose 
    \begin{equation*}
        \max_{\bm{y}\in\mc{I}(\mc{T})}g\big(\bm{y},\mc{T}\big) = \lambda.
    \end{equation*}
    In this case there exists $\bm{y}\in\mc{I}\big(\mc{T}\big)$ such that
    $J_{\lambda}\big(\mc{T}_{\bm{y}}\big) = J_{\lambda}(\bm{y})$.  We
    consider the cases when $\mc{T}_{\bm{y}\gamma}(\lambda)$ and
    $\mc{T}_{\bm{y}\gamma}$ are the same for all $\gamma\in\mc{Y}$ and
    when they differ for at least one $\gamma\in\mc{Y}$.  If
    $\mc{T}_{\bm{y}\gamma}(\lambda) = \mc{T}_{\bm{y}\gamma}$ for all
    $\gamma\in\mc{Y}$ then 
    \begin{equation*}
        \sum_{\gamma\in\mc{Y}}J_{\lambda}\big(\mc{T}_{\bm{y}\gamma}(\lambda)\big)
        = J_{\lambda}\big(\mc{T}_{\bm{y}}\big)
        = J_{\lambda}(\bm{y}),
    \end{equation*}
    and thus $\mc{T}(\lambda) \neq \mc{T}$ by
    Lemma~\ref{thm:optimal}.  If $\mc{T}_{\bm{y}\gamma}(\lambda) \neq
    \mc{T}_{\bm{y}\gamma}$ for at least one $\gamma\in\mc{Y}$ then
    $\mc{T}(\lambda) \neq \mc{T}$ again by Lemma~\ref{thm:optimal}.

    Finally, when 
    \begin{equation*}
        \max_{\bm{y}\in\mc{I}(\mc{T})}g\big(\bm{y},\mc{T}\big) > \lambda
    \end{equation*}
    then $\mc{T}(\lambda) \neq \mc{T}$ follows from the previous case and
    Lemma~\ref{thm:subtree}.
\end{proof}

\begin{figure}
    \begin{center}
        \input{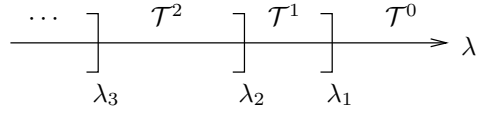}
    \end{center}
    \caption{ 
    \label{fig:exa2}
    For all $m \in \{0,1,\ldots, M-1\}$ the tree $\mc{T}^m$ is the
    smallest tree minimizing the cost $J_\lambda(\cdot)$ for any
    $\lambda\in (\lambda_{m+1},\lambda_m]$.}
\end{figure}

Let $\mc{T}^{0}$ denote the complete tree of depth $\kappa$.  Starting with
$\lambda_0=\infty$, for $m=\{1,\ldots, M\}$ recursively define
\begin{align*}
    \lambda_m 
    & \defeq \inf\{\lambda\leq  \lambda_{m-1}: \mc{T}^{m-1}(\lambda) = \mc{T}^{m-1}\:\}, \\
    \mc{T}^{m} 
    & \defeq \mc{T}^{m-1}(\lambda_m),
\end{align*}
where $M$ is the smallest integer such that $\lambda_{M+1}=0$, and with
$\lambda_1\defeq\infty$ if the set over which the infimum is taken is
empty. Lemma~\ref{thm:subtree} implies that for two consecutive
transition points $\lambda_{m}$ and $\lambda_{m+1}$, we have
$\mc{T}^0(\lambda)=\mc{T}^0(\lambda_m)$ for all $\lambda \in
(\lambda_{m+1},\lambda_m]$ as shown in Figure \ref{fig:exa2}. 

The following corollary is a consequence of Lemma \ref{thm:subtree} and Theorem
\ref{thm:lambda}. 
\begin{corollary}
    \label{thm:lambdacor} 
    For  $m\in\{1,\ldots,M\}$ 
    \begin{align}
        \label{eq:lambdacor1}
        \lambda_m 
        & = \max_{\bm{y}\in \mc{I}(\mc{T}^{m-1})}g(\bm{y},\mc{T}^{m-1}), \\
        \label{eq:lambdacor2}
        \mc{T}^{m} 
        & = \mc{T}^{m-1} \setminus 
        \bigcup_{\substack{\bm{y}\in \mc{I}(\mc{T}^{m-1}):\\
        g(\bm{y},\mc{T}^{m-1})=\lambda_m}}\mc{D}(\mc{T}^{m-1}_{\bm{y}}).
    \end{align}
    Moreover, the set $\{(\alpha_m,d_m)\}_{m=1}^M$ with
    \begin{align*}
        \alpha_m & \defeq \Pp(T(\mc{T}^m)<S), \\
        d_m & \defeq \E(T(\mc{T}^m)-S)^+,
    \end{align*}
    are the break-points of $d(\alpha)$.
\end{corollary}
\begin{proof}
    Let $\mc{T}^{m-1}$ be fixed. Equation~\eqref{eq:lambdacor1} follows
    directly from Theorem~\ref{thm:lambda}. For~\eqref{eq:lambdacor2},
    notice that as $J_\lambda(\mc{T})$ is continuous in $\lambda$, the
    definition of $\lambda_m$ yields
    $J_{\lambda_m}(\mc{T}^{m-1})=J_{\lambda_m}(\mc{T}^{m})$. Hence
    $\mc{T}^m$ is the smallest subtree of $\mc{T}^{m-1}$ with same root,
    and having a cost equal to $J_{\lambda_m}(\mc{T}^{m-1})$. From
    \eqref{eq:gequiv} and Lemma~\ref{thm:optimal}, we deduce that
    $\mc{T}^m$ is obtained from $\mc{T}^{m-1}$ by removing the
    descendants of any $\bm{y}\in \mc{I}(\mc{T}^{m-1})$ such that
    $g(\bm{y},\mc{T}^{m-1})=\lambda_m$.

    It remains to show that $\{(\alpha_m,d_m)\}_{m=1}^M$ are the
    break-points of $d(\alpha)$. By Lemma~\ref{thm:dual}, the
    break-points are achieved by non-randomized s.t.'s.  By Lemma
    \ref{thm:subtree} we have $\mc{T}^m=\mc{T}^0(\lambda_m)$, i.e.,
    $\mc{T}^m$ is the smallest subtree of $\mc{T}^0$ having the same
    root and minimizing the cost $J_{\lambda_m}(\mc{T})$. Hence, among
    the minimizers of $J_{\lambda_m}(\mc{T})$, $\mc{T}^m$ yields the
    largest $\Pp(T(\mc{T})<S)$. Therefore each pair $(\alpha_m,d_m)$ is
    a break-point. Conversely, given a break-point of $d(\alpha)$, let
    $\mc{T}$ be the smallest subtree of $\mc{T}^0$ achieving it.  Then
    $\mc{T}=\mc{T}^0(\lambda)$ for some $\lambda$.  Since
    $\mc{T}^0(\lambda_m)=\mc{T}^m$ we have that
    $\{\mc{T}^0(\lambda)\}_{\lambda \in
    {\mathbb{R}}}=\{\mc{T}^m\}_{m=0}^M$, and therefore $\mc{T}=\mc{T}^m$
    for some $m\in \{1,\ldots,M\}$.
\end{proof}

From Corollary \ref{thm:lambdacor},  we deduce the algorithm below that
fully characterizes $d(\alpha)$ by computing its set of break-points 
$\{(\alpha_m,d_m)\}_{m=1}^M$. 

\begin{center}
    \begin{minipage}{11cm}
        \rule{\textwidth}{0.5pt}
        Algorithm: Compute the break-points $\{(\alpha_m,d_m)\}_{m=1}^M$ of
        $d(\alpha)$
        \rule[0.15cm]{\textwidth}{0.5pt}
        \vspace{-0.6cm}
        \algsetup{indent=2em}
        \begin{algorithmic}
            \STATE $m \Leftarrow 0$
            \STATE $\lambda_0 \Leftarrow \infty$
            \STATE $\mc{T}^0 \Leftarrow \textrm{complete tree of depth $\kappa$}$
            \REPEAT
            \STATE $m \Leftarrow m+1$
            \STATE $\lambda_{m} \Leftarrow \max_{\bm{y}\in\mc{I}(\mc{T}^{m-1})}
            g\big(\bm{y},\mc{T}^{m-1}\big)$
            \STATE $\mc{T}^{m} \Leftarrow 
            \mc{T}^{m-1} \setminus \bigcup_{\substack{\bm{y}\in\mc{I}(\mc{T}^{m-1}):\\
            g(\bm{y},\mc{T}^{m-1})=\lambda_m}}\mc{D}(\mc{T}^{m-1}_{\bm{y}})$
            \STATE $\alpha_m \Leftarrow \Pp(T(\mc{T}^m)<S)$ 
            \STATE $d_m \Leftarrow \E(T(\mc{T}^m)-S)^+$
            \UNTIL{$\lambda_m=0$}
            \STATE $M \Leftarrow m-1$
        \end{algorithmic}
        \vspace{-0.2cm}
        \rule{\textwidth}{0.5pt}
    \end{minipage}
\end{center}

As a $\vert\mc{Y}\vert$-ary tree has less than
$\vert\mc{Y}\vert^{\kappa}$ non-terminal nodes, the algorithm terminates
after at most that many iterations. Further, one may check that each
iteration has a running time that is $\exp(O(\kappa))$.  Therefore, the
worst case running time of the algorithm is $\exp(O(\kappa))$. This is
to be compared, for instance, with exhaustive search that has a
$\Omega(\exp\exp(\kappa))$ running time (because all
break-points of $d(\alpha)$ are achieved by non-randomized s.t.'s and
there are already $2^{\vert\mc{Y}\vert^{\kappa-1}}$
$\vert\mc{Y}\vert$-ary trees having leaves at either depth $\kappa$ or
$\kappa-1$).  

In Sections~\ref{sec:perm} and~\ref{sec:lookahead} we will see that,
under certain conditions on $\{(X_i,Y_i)\}_{i\geq 1}$
and $S$, the running time of the algorithm is only
\emph{polynomial} in $\kappa$.

\subsection{A Lower Bound on the Reaction Delay}
\label{sec:bound}

From Corollary \ref{thm:lambdacor}, we may also deduce a lower bound on $d(\alpha)$.
Since $d(\alpha)$ is convex, we can lower bound it as 
\begin{equation}
    \label{eq:bound}
    d(\alpha) \geq d(0)+\alpha d'(0+)
\end{equation}
where $d '(0+)$ denotes the right derivative of $d $ at $\alpha=0$. By
Corollary~\ref{thm:lambdacor}, if $\lambda_1 < \infty$ then
$d (0)$ is achieved by the complete tree $\mc{T}^0$, and if $\lambda_1 =
\infty$ then $d (0)$ is achieved by $\mc{T}^1$ which is a strict subtree of
$\mc{T}^0$. Hence~\eqref{eq:bound} can be written as
\begin{equation}
    \label{eq:bound2}
    d (\alpha) \geq d (0)-
    \begin{cases}
        \alpha\lambda_1 & \text{if $\lambda_1<\infty$}, \\
        \alpha\lambda_2 & \text{else}.
    \end{cases}
\end{equation} 
Note that the above bound is tight for $\alpha\leq \alpha_1$ with $\alpha_1>0$ when
$\lambda_1<\infty$, and is tight for $\alpha\leq \alpha_2$ with $\alpha_2>0$ when
$\lambda_1=\infty$. The following example illustrates this bound.

\begin{example}
    \label{eg:bound}
    Let $\{X_i\}_{i\geq 1}$ be i.i.d.  $\textrm{Bernoulli}(1/2)$, and
    let the $Y_i$'s be the output of a binary symmetric channel with
    crossover probability $p\in (0,1/2)$ for input $X_i$.  Consider the
    s.t. $S$ defined as
    \begin{equation*}
        S \defeq 
        \begin{cases} 
            1 & \text{if $X_1 = 1$,} \\
            \kappa & \text{else.}
        \end{cases}
    \end{equation*}
    For $\kappa=2$, the tree corresponding to this s.t. 
    is depicted in Figure~\ref{fig:tree}.

    Since $p\in(0,1/2)$, it is clear that whenever $\mc{T}$ is not
    the complete tree of depth $\kappa$, we have $\Pp(T(\mc{T})<S)>0$,  hence
    \begin{equation*}
        d(0) = \E(T(\mc{T}^{0})-S)^+ = \frac{1}{2}(\kappa-1).
    \end{equation*}
    An easy computation using Corollary \ref{thm:lambdacor} yields
    \begin{equation*}
        \lambda_1 = \frac{1-p}{p}(\kappa-1),
    \end{equation*}
    and, using \eqref{eq:bound2}, we get
    \begin{equation}
        \label{eq:bound3}
        d(\alpha) \geq (\kappa-1)\Big(\frac{1}{2}-\alpha\frac{1-p}{p}\Big).
    \end{equation}
    Let us comment on \eqref{eq:bound3}. Consider any two correlated
    sequences $\{X_i\}_{i\geq 1}$ and $\{Y_i\}_{i\geq 1}$ and
    a s.t. $S$ with respect to the $X_i$'s. Intuition tells
    us that there are two factors affecting $d(\alpha)$. The
    first is the correlation between the $X_i$'s and $Y_i$'s, in the
    above example parameterized by $p$. The lower the correlation, the higher
    $d(\alpha)$ will be. The second factor is the ``variability" of
    $S$, and might be characterized by the difference in terms of depth
    among the leaves having large probability to be reached. In the
    above example the ``variability" might be captured by $\kappa-1$, since
    with probability $1/2$ a leaf of depth $1$ is reached, and with
    probability $1/2$ a leaf of depth $\kappa$ is attained.
\end{example}

\begin{example} \label{eg:continu}
    We consider one-bit message feedback communication when the forward
    and the feedback channels are binary erasure channels with erasure
    probabilities $\varepsilon$ and $p$, respectively. We refer the
    reader to Example~\ref{eg:communication} in Section~\ref{sec:intro}
    for the general problem setting.  We use the following transmission
    scheme (which is optimal in the case of noiseless feedback). The
    decoder keeps sending $0$ over the feedback channel until time $S$,
    the first time a non-erasure occurs or $\kappa$ time units have
    elapsed. From that point on, the decoder sends $1$. The encoder
    keeps sending the message bit it wants to deliver until time $T$ (a
    stopping time with respect to the output of the feedback channel).
    Ideally, we would like to choose $T=S$. This is possible if
    the feedback is noiseless, i.e., $p=0$. If $p>0$, we want to track the
    decoding time $S$ as closely as possible. The constant $\kappa$
    plays here the role of a ``time-out.'' In the following, we assume
    that $\varepsilon,p\in(0,1)$.

    Let us focus on $d(\alpha)$. One can show that
    $\lambda_1=\infty$ and therefore the bound~\eqref{eq:bound2} becomes
    $d(\alpha)\geq d(0)-\alpha\lambda_2$, where $\lambda_2 =
    \max_{\bm{y}\in\mc{I}(\mc{T}^{1})} g\big(\bm{y},\mc{T}^{1}\big)$
    from Corollary \ref{thm:lambdacor}. A somewhat involved computation
    yields
    \begin{equation}
        \label{eq:becbeta}
        d(\alpha) 
        \geq \left(\frac{p}{1-p}-\varepsilon^{1-\kappa}\alpha\right)(1+o(1))
    \end{equation}
    as $\kappa\rightarrow \infty$.

    The delay $d(\alpha)$ is interpreted as the time it takes the encoder
    to realize that the decoder has made a decision.
    Equation~\eqref{eq:becbeta} relates this delay to the channel parameters
    $\varepsilon$ and $p$, the probability $\alpha$ of stopping retransmission
    too early, and the value of the ``time-out'' $\kappa$. For the communication
    scheme considered here, there are two events  leading to decoding errors.
    The event $\{X_{\kappa}=0\}$, indicating that only erasures were received by
    the decoder until time $\kappa$, and the event $\{T< S\}$, indicating that
    the encoder stopped retransmission before the decoder received a non
    erasure. In both cases the decoder will make an error with probability
    $1/2$. Hence the overall probability of error $\Pp(\mc{E})$ can be bounded
    as
    \begin{equation*}
        \max\{\alpha,\varepsilon^{\kappa}\}\leq 2\Pp(\mc{E})
        \leq \alpha+\varepsilon^{\kappa}.
    \end{equation*}
    It is then reasonable to choose $\kappa=\frac{\log \alpha}{\log
    \varepsilon}$, i.e., to scale $\kappa$ with $\alpha$ so that both
    sources of errors have the same weight. This results in a delay of 
    \begin{equation*}
        d (\alpha) \geq \left(\frac{p}{1-p}-\varepsilon \right)(1+o(1))
    \end{equation*}
    as $\alpha \rightarrow 0$.

    Now suppose that the communication rate $R$ is computed with respect
    to the delay from the time communication starts until the time the
    decoder has made a decision and the encoder has realized this, i.e.,
    $\E S+\E (T-S)^+=\E(\max\{S,T\})$.  We conclude that the ``send
    until a non-erasure'' strategy asymptotically achieves a rate that
    is upper bounded as 
    \begin{equation*} 
        R  \leq 
        \frac{1}{\frac{1}{1-\varepsilon}+\frac{p}{1-p}-\varepsilon}\;. 
    \end{equation*}
    When $\varepsilon<p/(1-p)$, our bound is strictly below the capacity
    of the binary erasure channel $1-\varepsilon$. Hence $1/(1+\varepsilon)$
    represents a critical value for the erasure probability $p$ of the feedback
    channel above which the ``send until non-erasure'' strategy is strictly
    suboptimal. Indeed there exist block coding strategies, making no use of
    feedback, that (asymptotically) achieve rates up to $1-\varepsilon$, the
    capacity of the forward channel.  
\end{example}

\section{Permutation Invariant Stopping Times}
\label{sec:perm}

We consider a special class of s.t.'s and processes
$\{(X_i,Y_i)\}_{i\geq 1}$ for which the optimal tradeoff curve
$d(\alpha)$ and the associated optimal s.t.'s can be computed in
polynomial time in $\kappa$.

A s.t. $S$ with respect to $\{X_i\}_{i\geq 1}$ is
\emph{permutation invariant} if
\begin{equation*}
    \Pp(S\leq n\vert X^n=x^n) = \Pp(S\leq n\vert X^n=\pi(x^n))
\end{equation*}
for all permutations $\pi:\mc{X}^n\to\mc{X}^n$, all $x^n\in\mc{X}^n$ and
$n\in\{1,\ldots,\kappa\}$. Examples of permutation invariant s.t.'s are
$\inf\{i:X_i>c\}$ or $\inf\{i:\sum_{k=1}^i X_k >c\}$ for some constant
$c$ and assuming the $X_i$'s are positive.  The notion of a
permutation invariant s.t. is closely related to (and in fact slightly
stronger than) that of an exchangeable s.t. as defined in~\cite{LI}.

The following theorem establishes a key result, from which the running
time of one iteration of the algorithm can be deduced.

\begin{theorem}
    \label{thm:perm}
    Let $\{(X_i,Y_i)\}_{i\geq 1}$ be i.i.d. and $S$ be a
    permutation invariant s.t. with respect to $\{X_i\}_{i\geq 1}$.
    If $T(\mc{T})$ is non-randomized and permutation invariant then 
    \begin{equation*}
        g(\bm{y},\mc{T}) = g(\pi(\bm{y}),\mc{T})
    \end{equation*}
    for all $\bm{y}\in\mc{I}(\mc{T})$ and all permutations $\pi$.
\end{theorem}

We first establish two lemmas that will be used in the proof of
Theorem~\ref{thm:perm}. 

\begin{lemma}
    \label{thm:lemmatree}
    Let $T$ be a non-randomized s.t.
    with respect to $\{Y_i\}_{i\geq 1}$ and $\mc{T}$ the
    corresponding tree. Then $T$ is permutation invariant if and only if
    for all $\bm{y}\in\mc{I}(\mc{T})$ and permutations $\pi$,
    $\pi(\bm{y})\in\mc{I}(\mc{T})$.
\end{lemma}
\begin{proof}
    Assume $T$ is permutation invariant and let
    $y^n\in\mc{I}(\mc{T})$. Then 
    \begin{equation*}
        0 = \Pp(T\leq n \vert Y^n=y^n) 
        = \Pp(T\leq n\vert Y^n=\pi(y^n)),
    \end{equation*}
    and hence $\pi(y^n)\in\mc{I}(\mc{T})$.

    Conversely assume that, for all $\bm{y}\in\mc{I}(\mc{T})$ and
    permutations $\pi$, we have $\pi(\bm{y})\in\mc{I}(\mc{T})$. Pick an
    arbitrary $y^n$. First, if $\Pp(T\leq
    n\vert Y^n=y^n) = 0$, then $y^n\in\mc{I}(\mc{T})$, and by
    assumption also $\pi(y^n)\in\mc{I}(\mc{T})$. Thus
    $\Pp(T\leq n\vert Y^n=\pi(y^n))=0$. Second, if $\Pp(T\leq n\vert
    Y^n=y^n)=1$, then $y^n\notin\mc{I}(\mc{T})$, and by
    assumption also $\pi(y^n)\notin\mc{I}(\mc{T})$. Thus
    $\Pp(T\leq n\vert Y^n=\pi(y^n))=1$. 
\end{proof}

\begin{lemma}
    \label{thm:lemmay}
    Let $\{(X_i,Y_i)\}_{i\geq 1}$ be i.i.d. and $S$ be a
    permutation invariant s.t. with respect to
    $\{X_i\}_{i\geq 1}$. Then $S$ is a permutation invariant
    s.t. with respect to $\{Y_i\}_{i\geq 1}$.
\end{lemma}
\begin{proof}
    Using that the $\{(X_i,Y_i)\}_{i\geq 1}$ are i.i.d., 
    one can easily check that $S$ is a s.t. with respect to
    $\{Y_i\}_{i\geq 1}$. It remains to show that it is permutation
    invariant. For any permutation $\pi:\mc{X}^n\to\mc{X}^n$ 
    \begin{align*}
        \Pp(& S\leq n \vert Y^n=y^n) \\
        & = \sum_{x^n\in\mc{X}^n}\Pp(S\leq n\vert X^n=x^n)
        \Pp(X^n=x^n\vert Y^n=y^n) \\
        & = \sum_{x^n\in\mc{X}^n}\Pp(S\leq n\vert X^n=\pi^{-1}(x^n))\times \\
        & \quad \times\Pp(X^n=\pi^{-1}(x^n)\vert Y^n=y^n) \\
        & = \sum_{x^n\in\mc{X}^n}\Pp(S\leq n\vert X^n=x^n)
        \Pp(X^n=x^n\vert Y^n=\pi(y^n)) \\
        & = \Pp(S\leq n\vert Y^n=\pi(y^n)),
    \end{align*}
    where the second last equality follows by the permutation invariance of $S$
    and the fact that the $(X_i,Y_i)$'s are i.i.d. 
\end{proof}

\begin{proof}[Proof of Theorem~\ref{thm:perm}]
    We show that
    \begin{equation}
        \label{eq:gijo}
        g(\bm{y},\mc{T}) =
        \frac{b(\mc{T}_{\bm{y}})-b(\bm{y})}
        {a(\bm{y})-a(\mc{T}_{\bm{y}})}=g(\pi(\bm{y}),\mc{T})
    \end{equation}
    for all $\bm{y}\in\mc{I}(\mc{T})$. We prove that the numerator and the denominator in
    \eqref{eq:gijo} remain unchanged if we replace $\bm{y}$ by
    $\pi(\bm{y})$. Fix some $\bm{y}=y^n\in\mc{I}(\mc{T})$, and, to
    simplify notation, set $l=l(\bm{\gamma})$ until the end of this
    proof. For the denominator, using Lemma~\ref{thm:lemmay}
    we obtain 
    \begin{align}
        \label{eq:perm_step1}
        a & (\bm{y})- a(\mc{T}_{\bm{y}}) \nonumber\\
        & \defeq a(y^n)  -\sum_{y^n\bm{\gamma}\in\mc{L}(\mc{T}_{y^n})} a(y^n \bm{\gamma}) \nonumber\\
        & =  \Pp(Y^n=y^n)\Pp(S>n\vert Y^n=y^n) 
        \sum_{y^n\bm{\gamma}\in\mc{L}(\mc{T}_{y^n})}
        \Pp(Y^{n+l} = y^n\bm{\gamma})\Pp(S>n+l)\vert Y^{n+l}\!=y^n\bm{\gamma}) \nonumber \\
        & =  \Pp(Y^n=\pi(y^n))\Pp(S>n\vert Y^n=\pi(y^n))
        \sum_{y^n\bm{\gamma}\in\mc{L}(\mc{T}_{y^n})}
        \Pp(Y^{n+l} = \pi(y^n)\bm{\gamma})
        \Pp(S>n+l)\vert Y^{n+l} = \pi(y^n)\bm{\gamma}).
    \end{align}
    A consequence of Lemma~\ref{thm:lemmatree} is that the set of all
    $\bm{\gamma}$ such that $y^n\bm{\gamma}\in\mc{L}(\mc{T}_{y^n})$ is
    identical to the set of all $\bm{\gamma}$ such that
    $\pi(y^n)\bm{\gamma}\in\mc{L}(\mc{T}_{\pi(y^n)})$. Hence
    by~\eqref{eq:perm_step1}
    \begin{equation*}
        a(y^n)-a(\mc{T}_{y^n})
        = a(\pi(y^n))-a(\mc{T}_{\pi(y^n)}).
    \end{equation*}

    For the numerator in~\eqref{eq:gijo}, we have
    \begin{align}
        \label{eq:perm_step2}
        b(&\mc{T}_{y^n})-b(y^n) \nonumber\\
        & = \sum_{y^n\bm{\gamma}\in\mc{L}(\mc{T}_{y^n})}
        \Pp(Y^{n+l}=y^n\bm{\gamma})
        \Big(
        \E\big( (n+l-S)^+\big\vert Y^{n+l}=y^n\bm{\gamma} \big)
        -\E\big( (n-S)^+\big\vert Y^{n+l}=y^n\bm{\gamma} \big)
        \Big).
    \end{align}
    By Lemma~\ref{thm:lemmay}
    \begin{align*}
         \E\big( (n+l-S)^+ & \big\vert Y^{n+l}=y^n\bm{\gamma} \big)
        -\E\big( (n-S)^+\big\vert Y^{n+l}=y^n\bm{\gamma} \big)  \\
        & = \sum_{k=n}^{n+l-1}
        \Pp(S\leq k\vert Y^{n+l}=y^n\bm{\gamma} ) \\
        & = \sum_{k=n}^{n+l-1}
        \Pp(S\leq k\vert Y^{n+l}=\pi(y^n)\bm{\gamma} ).
    \end{align*}
    Combining this with~\eqref{eq:perm_step2} and using
    Lemma~\ref{thm:lemmatree} as before, we get
    \begin{equation*}
        b(\mc{T}_{y^n})-b(y^n) 
        = b(\mc{T}_{\pi(y^n)})-b(\pi(y^n)),
    \end{equation*}
    concluding the proof.
\end{proof}

We now show that one iteration of the algorithm has only polynomial
running time in $\kappa$. Specifically, we evaluate the running time to
compute $\mc{T}^{m+1}$ from $\mc{T}^m$ if $S$ and $T(\mc{T}^m)$ are
permutation invariant and if the $(X_i,Y_i)$'s are i.i.d. To that aim, we
assume the input of the algorithm to be in the form of a list of the
probabilities $\Pp(S\leq n\vert X^n=x^n)$ for all $x^n\in\mc{X}^n$ and
$n\in\{1,\ldots,\kappa\}$ --- specifying $S$ --- and a list of
$\Pp(X=x,Y=y)$ for all $x\in\mc{X}$ and $y\in\mc{Y}$ --- characterizing
the distribution of the process $\{(X_i,Y_i)\}_{i\geq 1}$. Note that as
$S$ is permutation invariant, we only have to specify $\Pp(S\leq n\vert
X^n=x^n)$ for each composition (or type) of $x^n$. Since the number of
compositions of length at most $\kappa$ is upper bounded by
$(\kappa+1)^{1+\vert\mc{X}\vert}$ --- any element $x\in\mc{X}$ appears
at most $k$ times in a string of length $k$ --- the list of these
probabilities has only polynomial size in $\kappa$.  Using a hash table,
we assume that, given $x^n$, the element $\Pp(S\leq n\vert X^n=x^n)$ in
the list can be accessed in $O(\kappa)$ time. The proof of the following
theorem is deferred to the appendix. 
\begin{theorem}
    \label{thm:permg}
    Let $\{(X_i,Y_i)\}_{i\geq 1}$ be i.i.d., let $S$ and $T(\mc{T}^m)$
    be permutation invariant s.t.'s with respect to $\{X_i\}_{i\geq 1}$
    and $\{Y_i\}_{i\geq 1}$ respectively, and let $\alpha_{m}=
    \Pp(T(\mc{T}^m)<S)$ and $d_m= \E(T(\mc{T}^m)-S)^+$ be
    given. Then $\mc{T}^{m+1}$, $\alpha_{m+1}$, and $d_{m+1}$ can be
    computed in polynomial time in $\kappa$.
\end{theorem}
As a corollary of Theorem \ref{thm:permg}, we obtain the worst case
running time for computing the set of break-points
$\{(\alpha_m,d_m)\}_{m=1}^M$ together with the associated optimal s.t.'s
$\{\mc{T}^m\}_{m=0}^M$.

\begin{corollary}
    \label{thm:permcor}
    Let $\{(X_i,Y_i)\}_{i\geq 1}$ be i.i.d. and $S$ be a
    permutation invariant s.t. with respect to $\{X_i\}_{i\geq 1}$. 
    If all $\{\mc{T}^m\}_{m=0}^{M}$ are permutation invariant,
    then the algorithm has a polynomial running time in $\kappa$.
\end{corollary}
\begin{proof}
    By Theorem~\ref{thm:permg} we only have to bound the number of
    iterations of the algorithm. To this end note that by Theorem
    \ref{thm:perm} every composition of $\bm{y}$ can be only once a
    maximizer of $g(\bm{y},\mc{T}^m)$ (as the corresponding nodes will
    be leaves in the next iteration of the algorithm). Hence, there are
    at most $O( (\kappa+1)^{1+\vert\mc{Y}\vert})$ iterations.  
\end{proof}

Note that, in the cases where $\{\mc{T}^m\}_{m=0}^{M}$ are not
permutation invariant, one may still be able to derive a lower bound on
$d(\alpha)$ in polynomial time in $\kappa$, using \eqref{eq:bound2}.
Indeed, the tree $\mc{T}^0$ is permutation invariant since it is
complete and, by Theorem~\ref{thm:permg}, if $\{(X_i,Y_i)\}_{i\geq 1}$
are i.i.d. and $S$ is permutation invariant, then the first subtree
$\mc{T}^1$ can be computed in polynomial time in $\kappa$. Therefore the
bound
\begin{equation}\label{eq:bound5}
    d(\alpha) \geq d(0)-\alpha\lambda_1
\end{equation}
can always be evaluated in polynomial time in $\kappa$ when the $(X_i,Y_i)$'s
are i.i.d. and $S$ is permutation invariant. Note that this bound is in
general weaker than the one derived in Section~\ref{sec:bound}. However,
when $\lambda_1<\infty$ the bound \eqref{eq:bound5} is tight for
$\alpha\in [0,\alpha_1]$ for some $\alpha_1>0$. It is easily checked that
the condition $\lambda_1<\infty$ is satisfied if $\Pp(S=\kappa,
Y^{\kappa-1}=y^{\kappa-1})>0$ for all $y^{\kappa-1}$.

In the next section, we present two examples for which the conditions of
Corollary~\ref{thm:permcor} are satisfied, and hence for which the
algorithm has a polynomial running time in $\kappa$. First, we consider
a TST problem that indeed can be formulated as a Bayesian change-point
problem.  Second, we consider the case of a pure TST problem, i.e.,
one that cannot be formulated as a Bayesian change-point problem. For
both examples, we provide an analytical solution of the Lagrange
minimization problem $\min_{T\leq \kappa} J_\lambda(T)$. 

\section{One-step Lookahead Stopping Times}
\label{sec:lookahead}

In this section, we show that under certain conditions the s.t.
that minimizes the Lagrangian $J_\lambda(T)$ can be found in closed form.

Define 
\begin{equation*}
    \mc{A}_n \defeq \Big\{ y^n\in\mc{Y}^n:
    \sum_{\gamma\in\mc{Y}}J_{\lambda}(y^n \gamma) \geq J_{\lambda}(y^n)
    \Big\},
\end{equation*}
and let
\begin{equation}
    \label{eq:lookahead}
    T^*_{\lambda} \defeq \min\big\{
    \kappa, \inf\{n: Y^n\in\mc{A}_n\}
    \big\}.
\end{equation}
In words, $T^*_{\lambda}$ stops whenever the current cost
\begin{equation*}
    \E((n-S)^+|Y^n=y^n)+\lambda \Pp(S>n|Y^n=y^n)
\end{equation*}
is less than the expected cost at time $n+1$, i.e.,
\begin{equation*}
    \E((n+1-S)^+|Y^n=y^n)+\lambda \Pp(S>n+1|Y^n=y^n)\;.
\end{equation*}
Recall that $\mc{T}^0$ denotes the complete tree of depth $\kappa$.  For
$(X_i,Y_i)$'s i.i.d., Theorem \ref{thm:lookahead} provides a sufficient
condition on $S$ for which $T(\mc{T}^0(\lambda))=T^*_{\lambda}$. In
words, the s.t.  $T^*_{\lambda}$ minimizes $J_\lambda(T)$ among all
s.t.'s bounded by $\kappa$. Furthermore, among all stopping times
minimizing $J_\lambda(T)$, the s.t. $T^*_{\lambda}$ admits the smallest
tree representation.  The proof of Theorem \ref{thm:lookahead} is
deferred to the appendix.

\begin{theorem}
    \label{thm:lookahead}
    Let $\{(X_i,Y_i)\}_{i\geq 1}$ be i.i.d., and let $S$ be a s.t. with
    respect to $\{X_i\}_{i\geq 1}$ that satisfies
    \begin{equation}\label{eq:monotone}
        \Pp(S=n\vert Y^{n-1})
        \geq \Pp(S=n+1\vert Y^{n})
    \end{equation}
    for all $n\in\{2,\ldots,\kappa\}$. Then 
    \begin{equation*}
        T(\mc{T}^0(\lambda)) = T^*_{\lambda}.
    \end{equation*} 
\end{theorem} 
Note that, unlike the algorithm, Theorem \ref{thm:lookahead} provides an
analytical solution only to the
inner minimization problem in \eqref{eq:dual}. To find the reaction delay
$d(\alpha)$ one still needs to maximize over the Lagrange multipliers $\lambda$.

Using Theorems~\ref{thm:permcor} and~\ref{thm:lookahead}, we now give
two examples of process $\{(X_i,Y_i)\}_{i\geq 1}$ and s.t. $S$ for which
the algorithm has only polynomial running time in $\kappa$.
\begin{example}
    \label{eg:lookahead1}
    Let $\{(X_i,Y_i)\}_{i\geq 1}$ be i.i.d. with the $X_i$'s taking
    values in $\{0,1\}$. Consider the s.t. $S\defeq\inf\{i: X_i=1\}$. We
    have for $n\geq 2$
    \begin{align*}
        \Pp(S=n \vert Y^{n-1})
        & =\Pp(S\geq n|Y^{n-1})\Pp(X_{n}=1) \\
        & \geq \Pp(S\geq n|Y^{n-1})\Pp(X_{n}=0|Y_{n}) \Pp(X_{n+1}=1) \\
        & = \Pp(S=n+1\vert Y^{n}).
    \end{align*}
    Hence, Theorem~\ref{thm:lookahead} yields that the one-step
    lookahead stopping time $T^*_{\lambda}$ defined in~\eqref{eq:lookahead}
    satisfies 
    \begin{equation*}
        T(\mc{T}^0(\lambda)) = T^*_{\lambda}.
    \end{equation*}

    We now show that the algorithm finds the set of break-points
    $\{(\alpha_m,d_m)\}_{m=0}^M$ and the corresponding
    $\{\mc{T}_m\}_{m=0}^M$ in polynomial running time in $\kappa$. To
    that aim, we first show that $T^*_{\lambda}$ is permutation
    invariant. By Lemma~\ref{thm:lemmatree}, we equivalently show that,
    for all $y^n$ and permutations $\pi$, if $y^n\notin\mc{A}_n$ then
    $\pi(y^n)\notin\mc{A}_n$. We have for $n<\kappa$
    \begin{align}
        \label{eq:looksym}
        \sum_{\gamma\in\mc{Y}} 
        J_{\lambda}(y^{n}\gamma)-J_{\lambda}(y^n)
        & = \Pp(Y^n=y^n)\Big( \Pp(S\leq n\vert Y^n=y^n)
        -\lambda\Pp(S=n+1\vert Y^n=y^n) \Big) \nonumber\\
        & = \Pp(Y^n=\pi(y^n))\Big( \Pp(S\leq n\vert Y^n=\pi(y^n))
        -\lambda\Pp(S=n+1\vert Y^n=\pi(y^n))
        \Big) \nonumber\\
        & = \sum_{\gamma\in\mc{Y}}J_{\lambda}(\pi(y^n)\gamma)-J_{\lambda}(\pi(y^n)),
    \end{align}
    where we have used Lemma~\ref{thm:lemmay} for the second equality.
    Thus $y^n\notin\mc{A}_n$ implies $\pi(y^n)\notin\mc{A}_n$,
    and therefore $T^*_{\lambda}$ is permutation invariant. Since
    $T(\mc{T}^0(\lambda)) = T^*_{\lambda}$ for all $\lambda\geq0$ by
    Theorem~\ref{thm:lookahead}, all $\{\mc{T}^m\}_{m=0}^{M}$ are
    permutation invariant.  Finally, because $S$ is permutation
    invariant, applying Corollary~\ref{thm:permcor} we conclude that the
    algorithm has indeed polynomial running time in $\kappa$.

    The problem considered in this example is actually a Bayesian
    change-point problem, as defined in Example~\ref{eg:cp} in Section
    \ref{sec:intro}. Here the change-point $\Theta\defeq S$ has
    distribution $\Pp(\Theta=n)=p(1-p)^{n-1}$, where $p\defeq\Pp(X=1)$.
    The conditional distribution of $Y_i$ given $\Theta$ is
    \begin{equation*}
        \Pp(Y_i=y_i\vert\Theta=n) =
        \begin{cases}
            \Pp(Y_i=y_i\vert X_i=0) & \text{if $i<n$}, \\
            \Pp(Y_i=y_i\vert X_i=1) & \text{if $i=n$}, \\
            \Pp(Y_i=y_i) & \text{if $i>n$}. \\
        \end{cases}
    \end{equation*}
    Note that, unlike the case considered by Shiryaev (see Example
    \ref{eg:cp} in Section \ref{sec:intro}), the distribution of the
    process at the change-point differs from the ones before and after
    it.
\end{example}

We now give an example that cannot be formulated as a change-point
problem and for which the one-step lookahead s.t. $T_\lambda^*$
minimizes the Lagrangian $J_\lambda(T)$.
\begin{example}
    \label{eg:lookahead2}
    Let $\{(X_i,Y_i)\}_{i\geq 1}$ be i.i.d. where the $X_i$'s and
    $Y_i$'s take values in $\{0,1\}$, and let $S\defeq \inf\{i\geq 1:
    \sum_{j=1}^iX_j=2\}$. A similar computation as for Example
    \ref{eg:lookahead1} reveals that if
    \begin{equation*}
        \Pp(X_i=1|Y_i)\geq \Pp(X_i=0|Y_i)
    \end{equation*}
    then Theorem \ref{thm:lookahead} applies, showing that the one-step
    lookahead stopping time $T^*_{\lambda}$ defined in~\eqref{eq:lookahead}
    satisfies $T(\mc{T}^0(\lambda)) = T^*_{\lambda}$. 

    Furthermore, since $S$ is permutation invariant, \eqref{eq:looksym}
    shows that $T^*_{\lambda}$ is permutation invariant. Applying
    Corollary \ref{thm:permcor}, one deduces that the algorithm has polynomial
    running time in $\kappa$ in this case as well. 

    The problem considered
    here is \emph{not} a change-point problem since, for $k>n$
    \begin{equation*}
        \Pp(S=k\vert Y^n=y^n,S > n) \neq \Pp(S=k\vert S>n),
    \end{equation*}
    and therefore~\eqref{eq:cpdiff} does not hold.
\end{example}

\section{Remarks}
\label{sec:conc}

In our study, we exploited the finite tree structure of bounded stopping times
defined over finite alphabet processes, and derived an algorithm that outputs
the minimum reaction delays for tracking a stopping time through noisy
observations, for any probability of false-alarm. This algorithm
has a complexity that is exponential in the bound of the stopping time we want
to track and, in certain cases, even polynomial. In comparison, an exhaustive
search has a complexity that is doubly exponential. 

The conditions under which the algorithm runs in polynomial time are,
unfortunately, not very explicit and require more study (see Corollary $10$). Explicit conditions,
however, are expected to be very restrictive on both the stochastic process and
the stopping time to be tracked. 

For certain applications, it is suitable to consider stopping times defined
over more general processes, such as continuous time over continuous alphabets.
In this case, how to solve the TST problem remains a wide open question. As a
first step, one might consider a time and alphabet quantization and apply our
result in order to derive an approximation algorithm.

\section*{Acknowledgments}
The authors are indebted to Professor G.~Wornell for providing them with
generous funding support without which this work would not have been possible.
We would also like to thank the anonymous reviewer for helpful comments.

\appendices

\section{Proof of Theorem \ref{thm:permg}}

In the following we write $\mc{T}$ for $\mc{T}^m$.  From
Theorem~\ref{thm:perm}, to find the $\bm{y}\in\mc{I}(\mc{T})$ maximizing
$g(\bm{y},\mc{T})$, we only have to compute $g(\bm{y},\mc{T})$ for all
possible compositions of $\bm{y}$.  The number of such compositions is
$O((\kappa+1)^{1+\vert\mc{Y}\vert})$.  We now show that
$g(\bm{y},\mc{T})$ can be computed in polynomial time in $\kappa$.
From the proof of Theorem~\ref{thm:perm}, we have to show that
$\Pp(S\leq n\vert Y^n=y^n)$ can be computed in polynomial time, and that
the sums in~\eqref{eq:perm_step1} and~\eqref{eq:perm_step2} can be
computed in polynomial time.

We have 
\begin{equation*}
    \Pp(S\leq n\vert Y^n=y^n)
    = \sum_{x^n\in\mc{X}^n}\Pp(S\leq n\vert X^n=x^n)\Pp(X^n=x^n\vert Y^n=y^n).
\end{equation*} 
Each term in the summation on the right hand side depends only on
the composition of $(x^n,y^n)$, and hence $\Pp(S\leq n\vert
Y^n=y^n)$ can be computed in polynomial time in $\kappa$.

Consider now the sum over all $y^n\bm{\gamma}\in\mc{L}(\mc{T}_{y^n})$
in~\eqref{eq:perm_step1}
\begin{equation}
    \label{eq:sumc}
    \sum_{y^n\bm{\gamma}\in\mc{L}(\mc{T}_{y^n})}
    a(y^n\bm{\gamma})
    = \sum_{y^n\tilde{\bm{\gamma}}\gamma\in\mc{L}(\mc{T}_{y^n})}
    a(y^n\tilde{\bm{\gamma}}\gamma).
\end{equation}
By Lemma~\ref{thm:lemmatree},
$y^n\tilde{\bm{\gamma}}\gamma\in\mc{L}(\mc{T}_{y^n})$ if and only if
$y^n\pi(\tilde{\bm{\gamma}})\gamma\in\mc{L}(\mc{T}_{y^n})$ for all
permutations $\pi$. And as
$a(y^n\tilde{\bm{\gamma}}\gamma)=a(y^n\pi(\tilde{\bm{\gamma}})\gamma)$,
we can compute~\eqref{eq:sumc} in polynomial time in $\kappa$.

Consider next the sum over all $y^n\bm{\gamma}\in\mc{L}(\mc{T}_{y^n})$
in~\eqref{eq:perm_step2}. Using Lemma~\ref{thm:lemmay}
\begin{multline*}
    \label{eq:sum2}
    \sum_{y^n\bm{\gamma}\in\mc{L}(\mc{T}_{y^n})}
    \sum_{k=n}^{n+l(\bm{\gamma})-1}
    \Pp(Y^{n+l(\bm{\gamma})}=y^n\bm{\gamma}) 
    \Pp(S\leq k\vert Y^{n+l(\bm{\gamma})}=y^n\bm{\gamma} ) \\
    = \sum_{y^n\bm{\gamma}\in\mc{I}(\mc{T}_{y^n})}
    \Pp(Y^{n+l(\bm{\gamma})}=y^n\bm{\gamma}) 
    \Pp(S\leq n+l(\bm{\gamma})\vert Y^{n+l(\bm{\gamma})}=y^n\bm{\gamma} ).
\end{multline*}
Applying Lemma~\ref{thm:lemmatree} as before, we conclude that the
right-hand side can be computed in polynomial time in $\kappa$.

It remains to prove that $\alpha_{m+1}$ and $d_{m+1}$ can be
computed in polynomial time in $\kappa$ from $\alpha_m$ and $d_m$.
This follows from the same argument, as it suffices to compute the
differences $b(\mc{T}_{\bm{y}^*})-b(\bm{y}^*)$ and
$a(\bm{y}^*)-a(\mc{T}_{\bm{y}^*})$ for all $\bm{y}^*$ maximizing
$g(\bm{y},\mc{T})$.

\section{Proof of Theorem \ref{thm:lookahead}}

Fix some $\lambda\geq 0$. Let us write $J_\lambda(T)$ as $\E(c(Y^T))$ where
\begin{equation*}
    c(y^n) \defeq 
    E( (n-S)^+\vert Y^n=y^n) + \lambda\Pp(S>n\vert Y^n=y^n).
\end{equation*}
We say that the $\{\mc{A}_n\}$ are \emph{nested} if, for any $n\geq 1$ and
$\gamma \in \mc{Y}$, we have that $y^n\in \mc{A}_n$ implies
$y^n\gamma\in \mc{A}_{n+1}$.  We show that \eqref{eq:monotone} implies
that the $\{\mc{A}_n\}$ are nested, and that this in turn implies that
the one-step lookahead stopping rule is optimal. The second part of the
proof is well known in the theory of optimal stopping and is referred as
the {\emph{monotone case}} (see, e.g., Chow et al.~\cite[Chapter
3]{CRS}). Here we provide an alternative proof that emphasizes the tree
structure of stopping times.

Note that $y^n\in\mc{A}_n$ if and only if $\E(c(Y^{n+1})\vert Y^n=y^n)\geq
c(y^n)$.  We now show that
\begin{equation}
    \label{eq:cost}
    \E(c(Y^{n+1})\vert Y^n)\geq c(Y^n) 
    \Longleftrightarrow \Pp(S\leq n\vert Y^n) \geq \lambda \Pp(S=n+1\vert Y^n).
\end{equation}
Since $\{(X_i,Y_i)\}_{i\geq 1}$ are i.i.d., $S$ is also a
(randomized) s.t. with respect to $\{Y_i\}_{i\geq 1}$ by Lemma \ref{thm:lemmay}.
It follows that
\begin{align*}
    c(Y^{n+1})
    & = E( (n+1-S)^+\vert Y^{n+1}) + \lambda\Pp(S>n+1\vert Y^{n+1}) \nonumber\\
    & = \sum_{k=1}^{n}\Pp(S\leq k\vert Y^{n+1})+\lambda\Pp(S>n+1\vert Y^{n+1}) \nonumber\\
    & = \sum_{k=1}^{n-1}\Pp(S\leq k\vert Y^{n})+\Pp(S\leq n\vert Y^n) \nonumber\\
    & \quad +\lambda\Pp(S>n\vert Y^{n}) -\lambda\Pp(S=n+1\vert Y^{n+1}) \nonumber\\
    & = c(Y^n) + \Pp(S\leq n\vert Y^n)-\lambda\Pp(S=n+1\vert Y^{n+1}),
\end{align*}
from which one deduces~\eqref{eq:cost}.

Next, we prove that the $\{\mc{A}_n\}$ are nested. By~\eqref{eq:cost} this is
equivalent to showing that, whenever for some $y^n$
\begin{equation}
    \label{eq:costhyp}
    \Pp(S\leq n\vert Y^n=y^n)
    \geq \lambda\Pp(S=n+1\vert Y^n=y^n),
\end{equation}
we also have 
\begin{equation} \label{eq:costhyp2}
    \Pp(S\leq n+1\vert Y^{n+1}=y^n\gamma)
    \geq \lambda\Pp(S=n+2\vert Y^{n+1}=y^n\gamma)
\end{equation}
for any $\gamma\in {\cal{Y}}$.  Suppose that~\eqref{eq:costhyp} holds
for some $y^n$. Using the fact that $S$ is a s.t.  with respect to the
$Y_i$'s (Lemma \ref{thm:lemmay}) together with the hypothesis of the
theorem yields for any $\gamma$
\begin{align*}
    \Pp (S\leq n+1 & \vert  Y^{n+1}=y^n\gamma)  - \lambda\Pp(S=n+2\vert Y^{n+1}=y^n\gamma) \\
    & \geq \Pp(S\leq n\vert Y^n=y^n) 
    - \lambda(S=n+2\vert Y^{n+1}=y^n\gamma) \\
    & \geq \lambda\Big(
    \Pp(S=n+1\vert Y^n=y^n)
    -\Pp(S=n+2\vert Y^{n+1}=y^n\gamma)
    \Big) \\
    & \geq 0,
\end{align*}
and therefore \eqref{eq:costhyp2} holds. Hence the $\{\mc{A}_n\}$ are nested.

Let $\mc{T}^*$ be the tree corresponding to $T^*_{\lambda}$.
The final step is to show that if the $\{\mc{A}_n\}$ are nested then
$\mc{T}^0(\lambda)=\mc{T}^*$. To that aim we show that
$\mc{I}(\mc{T}^*)\subset \mc{I}(\mc{T}^0(\lambda))$ and
$(\mc{I}(\mc{T}^*))^c\subset (\mc{I}(\mc{T}^0(\lambda)))^c$. Pick an
arbitrary $\bm{y}\in\mc{I}(\mc{T}^0)$. Using Lemma~\ref{thm:optimal}, we
compare $J_{\lambda}(\bm{y})$ with
$\sum_{\gamma}J_{\lambda}(\mc{T}^0_{\bm{y}\gamma}(\lambda))$.  We
distinguish two cases.  First suppose that $\bm{y}\in \mc{I}(\mc{T}^*)$,
i.e., $J_{\lambda}(\bm{y}) > \sum_{\gamma}J_{\lambda}(\bm{y}\gamma)$.
Then
\begin{equation*}
    J_{\lambda}(\bm{y})
    > \sum_{\gamma\in\mc{Y}}J_{\lambda}(\bm{y}\gamma)
    \geq \sum_{\gamma\in\mc{Y}}J_{\lambda}(\mc{T}^0_{\bm{y}\gamma}(\lambda)),
\end{equation*}
and hence $\bm{y}\notin\mc{L}(\mc{T}^0(\lambda))$. But since the
$\{\mc{A}_n\}$ are nested, no prefix of $\bm{y}$ can be an element of $\mc{L}(\mc{T}^0(\lambda))$ 
and hence $\bm{y}\in\mc{I}(\mc{T}^0(\lambda))$.

\balance

Second, assume $\bm{y}\notin \mc{I}(\mc{T}^*)$. If  $l(\bm{y})=\kappa$, then clearly
$\bm{y}\notin \mc{I}(\mc{T}^0(\lambda))$. If $l(\bm{y})<\kappa$,  then
$J_{\lambda}(\bm{y}) \leq \sum_{\gamma}J_{\lambda}(\bm{y}\gamma)$ and 
we now show by induction that this implies that
$\mc{T}^0_{\bm{y}}(\lambda)=\{\bm{y}\}$. Note first that as the
$\{\mc{A}_n\}$ are nested, we have for any
$\tilde{\bm{y}}\in\mc{I}(\mc{T}^0_{\bm{y}})$ (i.e., for any $\tilde{\bm{y}} $ with prefix $\bm{y}$)
\begin{equation}
    \label{eq:costind}
    J_{\lambda}(\tilde{\bm{y}}) 
    \leq \sum_{\gamma\in\mc{Y}}J_{\lambda}(\tilde{\bm{y}}\gamma).
\end{equation}
Assume first that $\mc{T}^0_{\tilde{\bm{y}}}$ has depth one.
Then~\eqref{eq:costind} implies by Lemma~\ref{thm:optimal} that
$\mc{T}^0_{\tilde{\bm{y}}}(\lambda) = \{\tilde{\bm{y}}\}$.  Suppose
then that this is true for all $\mc{T}^0_{\tilde{\bm{y}}}$ of depth at
most $k-1$. Let $\mc{T}^0_{\tilde{\bm{y}}}$ have depth $k$. Then by
the induction hypothesis and \eqref{eq:costind}
\begin{equation*}
    \sum_{\gamma\in\mc{Y}}J_{\lambda}(\mc{T}^0_{\tilde{\bm{y}}\gamma}(\lambda))
    = \sum_{\gamma\in\mc{Y}}J_{\lambda}(\tilde{\bm{y}}\gamma)
    \geq J_{\lambda}(\tilde{\bm{y}}),
\end{equation*}
and thus $\mc{T}^0_{\tilde{\bm{y}}}(\lambda)=\{\tilde{\bm{y}}\}$
by Lemma~\ref{thm:optimal}, concluding the induction step. 
This implies $\bm{y}\notin\mc{I}(\mc{T}^0(\lambda))$.

\bibliography{stopping_journal}

\end{document}